\documentclass{compositio}

\usepackage{amssymb, amsmath, rotate, multirow,epsfig}
\usepackage[all,2cell,dvips]{xy} \UseAllTwocells \SilentMatrices
\usepackage{graphicx}
\usepackage{color}

\newcommand\luijklist{\begin{itemize}
                      \setlength{\itemsep}{-1mm}
                     }
\newcommand\eindluijklist{\end{itemize}}
\newcommand\luijkreflist{\begin{itemize}
                      \setlength{\itemsep}{-1mm}
                     }
\newcommand\eindluijkreflist{\end{itemize}}

\renewcommand\check[1]{}
\newcommand\donecheck[1]{}

\newcommand\checkgone[1]{}

\makeatletter
\def\eqalign#1{\null\,\vcenter{\openup\jot\m@th
  \ialign{\strut\hfil$\displaystyle{##}$&$\displaystyle{{}##}$\hfil
      \crcr#1\crcr}}\,}
\makeatother

\newcommand\myref[1]{\cite{#1}}

\newtheorem{theorem}{Theorem}[section]
\newtheorem{proposition}[theorem]{Proposition}
\newtheorem{lemma}[theorem]{Lemma}
\newtheorem{corollary}[theorem]{Corollary}

\theoremstyle{definition}

\theoremstyle{remark}

\newtheorem{remark}{Remark}
\newtheorem{question}{Question}

\newcommand\lcm{\mathop{\rm lcm} \nolimits}
\newcommand\oX{\overline{X}}
\newcommand\oQ{\overline{\mathbb{Q}}}
\newcommand\ok{\overline{k}}
\newcommand\oZ{\overline{Z}}
\newcommand\oF{\overline{\F}}
\newcommand\PGL{\mathop{\rm PGL} \nolimits}

\renewcommand\O{\mathcal{O}}

\newcommand\Z{\mathbb{Z}}
\newcommand\Q{\mathbb{Q}}
\newcommand\F{\mathbb{F}}

\newcommand\C{\mathbb{C}}

\renewcommand\P{\mathbb{P}}

\newcommand\Spec{\mathop{\rm Spec} \nolimits}

\newcommand\Br{\mathop{\rm Br} \nolimits}

\newcommand\Tr{\mathop{\rm Tr} \nolimits}
\newcommand\disc{\mathop{\rm disc} \nolimits}

\newcommand\Pic{\mathop{\rm Pic} \nolimits}
\renewcommand\deg{\mathop{\rm deg} \nolimits}
\newcommand\NS{\mathop{\rm NS} \nolimits}
\newcommand\Hom{\mathop{\rm Hom} \nolimits}

\newcommand\tors{{\mathop{\rm tors} \nolimits}}

\newcommand\rk{\mathop{\rm rk} \nolimits}
\newcommand\Id{\mathop{\rm Id} \nolimits}

\newcommand\et{\text{\rm \'et}}

\newcommand\isom{\cong}

\renewcommand\tilde{\widetilde}

\newcommand\X{\mathfrak{X}}

\newcommand\myrefart[6]{\bibitem[#1]{#1}{#2, \emph{#3}, #4, #5, pp. #6.}}
\newcommand\myrefartfive[5]{\bibitem[#1]{#1}{#2, \emph{#3}, #4, #5.}}
\newcommand\myrefbook[5]{\bibitem[#1]{#1}{#2, \emph{#3} (#4, #5).}}

\newcommand\morimukai{MM83}
\newcommand\beau{Be85}
\newcommand\pjatshaf{PS71}
\newcommand\kloos{Kl05}
\newcommand\geem{VG04}
\newcommand\ellen{El04}
\newcommand\teras{Te85}
\newcommand\BT{BT00}
\newcommand\BPV{BPV84}
\newcommand\silv{Si86}
\newcommand\shiodaone{Sh81}
\newcommand\hag{Ha77}
\newcommand\nyog{NO85}
\newcommand\tatethree{Ta65}
\newcommand\maz{Ma77}
\newcommand\sgasix{Gr71}
\newcommand\milneone{Mi75}
\newcommand\milne{Mi80}
\newcommand\sgaseventwo{DK73}
\newcommand\borwein{Bo95}

\newcommand\luijkheron{VL04}

\begin{document}

\title[K3 surfaces with Picard number one]
{K3 surfaces with Picard number one and infinitely many rational points}
\author{Ronald van Luijk}
\email{rmluijk@math.berkeley.edu}
\address{University of California\\Berkeley, CA\\USA}
%

\begin{abstract}
In general, not much is known about the arithmetic of K3
surfaces. Once the geometric Picard number, which is the rank of the 
N\'eron-Severi group over an algebraic closure of the base field, is
high enough, more structure is known and more can be said.  
However, until recently not a single K3 surface was known to have
geometric Picard number one. We give
explicit examples of such surfaces over the rational numbers.
This solves an old problem that has been attributed to Mumford.
The examples we give also contain infinitely many rational points,
thereby answering a question of Swinnerton-Dyer and Poonen. 
\end{abstract}

\maketitle


\section{Introduction}
\label{sec:introduction}

K3 surfaces are the two-dimensional analogues of elliptic curves in
the sense that their canonical sheaf is trivial. However, as opposed
to elliptic curves, little is known about the arithmetic of K3
surfaces in general. It is for instance an open question if there
exists a K3 surface $X$ over a number field such that the set
of rational points on $X$ is neither empty, nor dense. We will answer
a longstanding question regarding the Picard group of a K3
surface. The Picard group of a K3 surface $X$ over a field $k$ is a
finitely generated free abelian group, the rank of which is called the Picard
number of $X$. The Picard number of $\oX=X \times_k \ok$, where $\ok$
denotes an algebraic closure of $k$, is called the
geometric Picard number of $X$. We will give the first known examples of
explicit K3 surfaces shown to have geometric Picard number $1$.

Bogomolov and Tschinkel \myref{\BT}
showed an interesting relation between the geometric Picard number
of a K3 surface $X$ over a number field $K$ and the arithmetic of $X$.
They proved that if the geometric Picard number is at least $2$, 
then in most cases the rational points on $X$ are potentially
dense, which means that there exists a finite field extension $L$ of
$K$ such that the set $X(L)$ of $L$-rational points is Zariski dense
in $X$, see \myref{\BT}. 
However, it is not yet known whether there exists any K3 surface
over a number field and with geometric Picard number $1$ on which the
rational points are potentially dense. Neither do we know if there exists
a K3 surface over a number field and with geometric Picard number $1$ on
which the rational points are {\em not} potentially dense! 

In December 2002, at the AIM workshop on rational and integral points on
higher-dimensional varieties in Palo Alto, Swinnerton-Dyer and 
Poonen asked a related question. \check{first names?}
They asked whether there exists a K3 surface over a number field and
with Picard number $1$ that contains infinitely many  
rational points. In this article we will show that such K3 surfaces do
indeed exist. It follows from our main theorem. 

\begin{theorem}\label{picmain}
In the moduli space of K3 surfaces polarized by a very ample divisor
of degree $4$, the set of  
surfaces defined over $\Q$ with geometric Picard number $1$ and
infinitely many rational points is Zariski dense. 
\end{theorem}

Note that
a polarization of a K3 surface is a choice of an ample divisor $H$.
The degree of such a polarization is $H^2$. A K3 surface
polarized by a very ample divisor of degree $4$ is a smooth quartic 
surface in $\P^3$. 
We will prove the main theorem by exhibiting an explicit family of quartic 
surfaces in $\P^3_\Q$ with geometric Picard number $1$ and
infinitely many rational points. 
Proving that these surfaces contain infinitely many rational
points is the easy part. It is much harder to prove that the geometric
Picard number of these surfaces equals $1$. 
It has been known since Noether that a general hypersurface in
$\P^3_\C$ of degree at least $4$ has geometric Picard number $1$. A
modern proof of this fact was given by Deligne, see 
\myref{\sgaseventwo},  Thm. XIX.1.2. Despite this fact, it has been an
old challenge, attributed to Mumford and disposed of in this article,
to find even one 
explicit quartic surface, defined over a number field, of which the
geometric Picard  
number equals $1$. Deligne's result does not actually imply that such
surfaces exist, as ``general'' means ``up to a countable union of
closed subsets of the moduli space.'' A priori, this could exclude all
surfaces defined over $\overline{\Q}$. Terasoma and 
Ellenberg have proven independently that such surfaces do exist. The
following theorems state their results.

\begin{theorem}[(Terasoma, 1985)]\label{tera}
For any positive integers $(n; a_1,\ldots,a_d)$ not equal to 
$(2;3), (n;2)$, or $(n;2,2)$, and with $n$ even, there is a smooth
complete intersection $X$ over $\Q$ 
of dimension $n$ defined by equations of degrees $a_1,\ldots,a_d$
such that the middle geometric Picard number of $X$ is $1$.
\end{theorem}
\begin{proof}
See \myref{\teras}.
\end{proof}

\begin{theorem}[(Ellenberg, 2004)]\label{ellen}
For every even integer $d$ there exists a number field $K$ and a
polarized K3 surface $X/K$ of degree $d$, with geometric Picard
number $1$.
\end{theorem}
\begin{proof}
See \myref{\ellen}.
\end{proof}

The proofs of Terasoma and Ellenberg are
ineffective in the sense that they do not give explicit examples. 
In principle it might be possible to extend their methods to test whether
a given explicit K3 surface has geometric Picard number $1$. In practice
however, it is an understatement to say that the amount of work
involved is not encouraging. The explicit examples we will give to
prove the main theorem also prove the case $(n; a_1,\ldots,a_d)=(2;4)$ of
Theorem \ref{tera} and the case $d=4$ of Theorem \ref{ellen}.

Shioda did find explicit examples of surfaces with geometric
Picard number $1$. In fact, he has shown that for every prime $m\geq
5$ the surface in $\P^3$ given by 
$$
w^m+xy^{m-1}+yz^{m-1}+zx^{m-1} = 0
$$
has geometric Picard number $1$, see \myref{\shiodaone}. However, for
$m=4$ this equation determines a K3 surface with maximal geometric Picard
number $20$, i.e., a singular K3 surface.

Before we prove the main theorem in Section \ref{sec:maintheorem}, we
will recall some definitions and results.

\section{Prerequisites}
\label{sec:prerequisites}

A {\em lattice} is a free $\Z$-module $L$ of finite rank, endowed with a
symmetric, bilinear, nondegenerate map $\langle
\underline{\,\,\,\,}\,,\underline{\,\,\,\,} \rangle \colon L \times L
\rightarrow \Q$, called the {\em pairing} of the lattice.
A {\em sublattice} of $L$ is a submodule $L'$ of $L$, such that the induced
bilinear pairing on $L'$ is nondegenerate. 
The {\em Gram matrix} of a lattice $L$ with respect to a given basis
$x=(x_1,\ldots, x_n)$ is $I_x = (\langle x_i,x_j \rangle)_{i,j}$.
The {\em discriminant} of $L$ is defined by $\disc L = \det I_x$ 
for any basis $x$ of $L$. For any sublattice $L'$ of finite 
index in $L$ we have $\disc L' = [L:L']^2 \disc L$. The image of $\disc L$ 
and $\disc L'$ in $\Q^*/{\Q^{*}}^2$ is the discriminant of the inner
product space $L_\Q$, where the inner product is induced by the
pairing of $L$.

Let $X$ be a smooth, projective, geometrically integral
surface over a field $k$ and set $\oX= X \times_k \ok$, where $\ok$
denotes an algebraic closure of $k$. 
The {\em Picard group} $\Pic X$ of $X$ is the group of line bundles on $X$
up to isomorphism, or equivalently, the group of divisor classes
modulo linear equivalence. The divisor classes that become algebraically
equivalent to $0$ over $\overline{k}$ (see \myref{\hag}, exc. V.1.7) 
form a subgroup $\Pic^0 X$ of $\Pic X$. The quotient
is the {\em N\'eron-Severi group} $\NS(X)=\Pic X/\Pic^0 X$, which is a
finitely generated abelian group, see \myref{\hag}, exc.~V.1.7--8,
or \myref{\milne}, Thm. V.3.25, for surfaces or \myref{\sgasix},
Exp. XIII, Thm. 5.1 in general.
The intersection pairing endows the group $\NS(X)/\NS(X)_\tors$ with
the structure of a lattice. Its rank is called the {\em Picard number}
of $X$. The Picard number of $\oX$ is called the {\em geometric Picard
  number} of $X$. 

By definition a smooth, projective, geometrically integral surface $X$
is a {\em K3 surface} if the canonical sheaf $\omega_X$ on $X$ is 
trivial and $H^1(\oX, \O_{\oX})=0$. Examples of K3 surfaces are smooth
quartic surfaces in $\P^3$. The Betti numbers of a K3 surface are
$b_0=1$, $b_1=0$, $b_2=22$, $b_3=0$, and $b_4=1$.

\begin{lemma}\label{NSeven}
If $X$ is a K3 surface, then $\Pic^0 X$ is trivial,  
the N\'eron-Severi group $\NS(X) \isom \Pic X$ is torsion free, and 
the intersection pairing on $\NS(X)$ is even.
\end{lemma}
\begin{proof}
See \myref{\BPV}, p. 21 and Prop. VIII.3.2. 
\end{proof}

For any scheme $Z$ over $\F_q$ with $q=p^r$ and $p$ prime and any
prime $l \neq p$, we define 
$$H^2_{\et}(Z,\Q_l)=\left(\lim_{\leftarrow} H^2_\et(Z,\Z/l^n\Z)\right)
\otimes_{\Z_l} \Q_l,$$ 
see \myref{\tatethree}, p. 94. 
Furthermore, for every integer $m$ and every vector space $H$ over
$\Q_l$ with the Galois group $G(\overline{\F_q}/\F_q)$ 
acting on it, we define the twistings of $H$
to be the $G(\overline{\F_q}/\F_q)$-spaces $H(m)=H \otimes_{\Q_l}
W^{\otimes m}$, where 
$$
W=\Q_l \otimes_{\Z_l} (\lim_{\leftarrow} \mu_{l^n})
$$
is the one-dimensional $l$-adic vector space on which
$G(\overline{\F_q}/\F_q)$ operates according to its action on the group 
$\mu_{l^n}\subset \overline{\F_q}$ of $l^n$-th roots of unity. Here we use 
$W^{\otimes 0} = \Q_l$ and
$W^{\otimes m} = \Hom(W^{\otimes -m}, \Q_l)$ for $m<0$.
For a surface $Z$ over $\overline{\F_q}$ the cup-product gives 
$H^2_{\et}(Z,\Q_l)(m)$ the structure of an inner product space for all
integers $m$.

Proposition \ref{nerred} describes the behavior of the
N\'eron-Severi group under good reduction. Its corollary will be used 
to show that the geometric Picard number of a certain surface
is equal to $1$.

\begin{proposition}\label{nerred}
Let $A$ be a discrete valuation ring of
a number field $L$ with residue field $k\isom \F_{q}$. Let $S$ be an integral
scheme with a morphism $S \rightarrow \Spec A$ that is projective and
smooth of relative dimension $2$. Assume that the
surfaces $\overline{S} = S_{\overline{L}}$ and $\tilde{S} =
S_{\overline{k}}$ are integral. Let $l \nmid q$ be a prime number.
Then there are natural injective homomorphisms
\begin{equation}\label{nerredmaps}
\NS(\overline{S}) \otimes \Q_l \hookrightarrow
\NS(\tilde{S}) \otimes \Q_l \hookrightarrow
H^2_{\et}(\tilde{S},\Q_l)(1)
\end{equation}
of finite dimensional inner product spaces over $\Q_l$. The second injection
respects the Galois action of $G(\overline{k}/k)$.
\end{proposition}
\begin{proof}
See \myref{\luijkheron}, Proposition 6.2.
\end{proof}

Recall that for any scheme $Z$ over $\F_q$ with $q=p^r$ and $p$ prime, 
the {\em absolute Frobenius} $F_Z\colon Z \rightarrow Z$ of $Z$ acts
as the identity on points, and by $f \mapsto f^p$ on the structure
sheaf. Set $\Phi_Z = F_{Z}^r$ and $\oZ = Z \times \oF_q$.
Let $\Phi^*_Z$ denote the automorphism on $H^2_{\et}(\oZ,\Q_l)$ induced by
$\Phi_Z \times 1$ acting on $Z \times \oF_q=\oZ$. 

\begin{corollary}\label{rankNSbound}
With the notation as in Proposition \ref{nerred},
the ranks of $\NS(\tilde{S})$ and
$\NS(\overline{S})$ are bounded from above by the number of
eigenvalues $\lambda$ of $\Phi^*_{S_{k}}$ for which 
$\lambda/q$ is a root of unity, counted with multiplicity.
\end{corollary}
\begin{proof}
By Proposition \ref{nerred} any upper bound for the rank of
$\NS(\tilde{S})$ is an upper bound for the rank of $\NS(\overline{S})$.
Let $\sigma$ denote the $q$-th power Frobenius map, i.e., the
canonical topological 
generator of $G(\overline{k}/k)$. For any positive integer $m$,
let $\sigma^*$ and $\sigma^*(m)$ denote the automorphisms induced on 
$\NS(\tilde{S})\otimes \Q_l$ and $H^2_\et(\tilde{S},\Q_l)(m)$ 
respectively. As all divisor classes are defined over some finite
extension of $k$, some power of Frobenius acts as the identity on 
$\NS(\tilde{S})$, so all eigenvalues of $\sigma^*$ acting on 
$\NS(\tilde{S})$ are roots of unity.
It follows from Proposition \ref{nerred} that
the rank of $\NS(\tilde{S})$ is bounded from above by the number of
roots of $\sigma^*(1)$ that are a root of unity. As the eigenvalues of
$\sigma^*(0)$ differ from those of $\sigma^*(1)$ by a factor of 
$q$, this equals the number of roots $\lambda$ of $\sigma^*(0)$ 
for which $\lambda q$ is a root of unity. The Corollary follows from
the fact that $\Phi^*_{S_{k}}$ acts on $H^2_{\et}(\oZ,\Q_l)$ as the
inverse of $\sigma^*(0)$. See also \myref{\luijkheron}, Corollary 6.3.
\end{proof}

\begin{remark}
Tate's conjecture states that the upper bound mentioned is
actually equal to the rank of $\NS(\tilde{S})$, see \myref{\tatethree}.
Tate's conjecture has been proven for ordinary K3 surfaces over fields
of characteristic $p \geq 5$, see \myref{\nyog}, Thm. 0.2.
\end{remark}

To find the characteristic polynomial of Frobenius as in Corollary
\ref{rankNSbound}, we will use the following lemma.

\begin{lemma}\label{piccharpoly}
Let $V$ be a vector space of dimension $n$ and $T$ a linear operator 
on $V$. Let $t_i$ denote the trace of $T^i$. Then the characteristic 
polynomial of $T$ is equal to 
$$
f_T(x) = \det(x \cdot \Id - T) = x^n+c_1x^{n-1}+c_2 x^{n-2}+\ldots +c_n,
$$
with the $c_i$ given recursively by 
$$
c_1 = -t_1 \qquad \mbox{and} \qquad -kc_k = t_k + \sum_{i=1}^{k-1}c_it_{k-i}.
$$
\end{lemma}
\begin{proof}
This is Newton's identity, see \myref{\borwein}, p. 5.
\end{proof}

\section{Proof of the main theorem}
\label{sec:maintheorem}

First we will give a family of smooth quartic surfaces 
in $\P^3$ with Picard number $1$. Let $R = \Z[x,y,z,w]$ be the
homogeneous coordinate ring of $\P^3_\Z$. Throughout the rest of this
article, for any homogeneous polynomial $h \in R$ of degree $4$, let
$\X_h$ denote the scheme in $\P^3_\Z$ given by  
\begin{equation}\label{pictheeq}
wf_1+2zf_2 = 3g_1g_2 +6h, 
\end{equation}
with $f_1,f_2,g_1,g_2\in R$ equal to 
$$
\eqalign{
f_1 \,=&\,\, x^3-x^2y-x^2z+x^2w-xy^2-xyz+2xyw+xz^2+2xzw+y^3+\cr
     &\,\,+y^2z-y^2w+yz^2+yzw-yw^2+z^2w+zw^2+2w^3,\cr
f_2 \,=&\,\, xy^2+xyz-xz^2-yz^2+z^3,\cr
g_1 \,=&\,\, z^2+xy +yz,\cr
g_2 \,=&\,\, z^2+xy. \cr
}
$$
Its base extensions to $\Q$ and $\oQ$ are
denoted $X_h$ and $\oX_h$ respectively.

\begin{theorem}\label{picexampletheorem}
Let $h \in R$ be a homogeneous polynomial of degree $4$. Then the quartic
surface $X_h$ is smooth over $\Q$ and has  
geometric Picard number $1$. The Picard group $\Pic \overline{X}_h$
is generated by a hyperplane section. 
\end{theorem}
\begin{proof}
For $p=2,3$, let $X_p/\F_p$ denote the fiber of $\X_h \rightarrow
\Spec \Z$ over $p$. As they are independent of $h$, one easily checks that
$X_p$ is smooth over $\F_p$ for $p=2,3$. 
As the morphism $\X_h \rightarrow \Spec \Z$ is flat and projective, 
it follows that the generic fiber $X_h$ of 
$\X_h \rightarrow \Spec \Z$ is smooth over $\Q$ as well, cf. \myref{\hag},
exc. III.10.2.  

We will first show that $X_2$ and $X_3$ have geometric Picard number $2$.
For $p=2,3$, let $\Phi_p$ denote the absolute Frobenius of $X_p$.
Set $\oX_p = X_p \times \oF_p$ and 
let $\Phi_p^*(i)$ denote the automorphism on $H^i_{\et}(\oX_p,\Q_l)$ induced by
$\Phi_p \times 1$ acting on $\oX_p = X_p \times_{\F_p} \oF_p$. Then by 
Corollary \ref{rankNSbound}
the geometric Picard number of $X_p$ is bounded 
from above by the number of eigenvalues $\lambda$ of
$\Phi^*_p(2)$ for which $\lambda/p$ is a root of unity. We will find
the characteristic polynomial of $\Phi_p^*(2)$ from the traces of its
powers. These traces we will compute with the Lefschetz formula
\begin{equation}\label{piclefschetz}
\# X_p(\F_{p^n}) = \sum_{i=0}^4 (-1)^i \Tr(\Phi_p^*(i)^n).
\end{equation}
As $X_p$ is a smooth hypersurface in $\P^3$ of degree $4$, 
it is a K3 surface and its Betti numbers are 
$b_0=1$, $b_1=0$, $b_2=22$, $b_3=0$, and $b_4=1$. It follows that 
$\Tr(\Phi_p^*(i)^n)=0$ for $i=1,3$, and for $i=0$ and $i=4$ the
automorphism $\Phi_p^*(i)^n$ has
only one eigenvalue, which by the Weil conjectures equals $1$ and
$p^{2n}$ respectively. From the Lefschetz formula (\ref{piclefschetz}) we
conclude $\Tr(\Phi_p^*(2)^n) = \# X_p(\F_{p^n})-p^{2n}-1$. After
counting points on $X_p$ over $\F_{p^n}$ for $n=1,\ldots,11$, this
allows us to compute the traces of the first $11$ powers of
$\Phi_p^*(2)$. With Lemma \ref{piccharpoly} we can then compute the first 
coefficients of the characteristic polynomial $f_p$ of $\Phi_p^*(2)$,
which has degree $b_2=22$. 
Writing $f_p = x^{22} + c_1 x^{21} + \ldots + c_{22}$ we find the
following table.

$$
\begin{array}{c||c|c|c|c|c|c|c|c|c|c|c}
p & c_1 & c_2 &c_3&c_4&c_5&c_6&c_7&c_8&c_9&c_{10}&c_{11} \cr
\hline
2 & -3  &  -2 & 12 & 0 &-32&64&-128 &128 & 256 &0 &-2048 \cr
\hline
3 &-5&-6&72 & 27 &-891 &0 &9477&-4374& -78732 & 19683 & 708588 \cr
\end{array}
$$
\smallskip

The Weil conjectures give a functional equation 
$p^{22}f_p(x) = \pm x^{22} f_p(p^2/x)$. As in our case 
(both for $p=2$ and $p=3$) the middle coefficient $c_{11}$ of $f_p$ is
nonzero, the sign of the functional equation is positive. This
functional equation allows us to compute the remaining coefficients of
$f_p$. 

If $\lambda$ is a root of $f_p$ then $\lambda/p$ is a root of
$\tilde{f}_p(x) = p^{-22}f_p(px)$. Hence, the number of roots of
$\tilde{f}_p(x)$ that are also a root of unity gives an upper bound 
for the geometric Picard number of $X_p$. After factorization into
irreducible factors, we find 
{ 
$$
\eqalign{
\tilde{f}_2 =& \textstyle{\frac{1}{2}}(x-1)^2\left(
2x^{20} + x^{19} - x^{18} + x^{16} + x^{14} + x^{11} + 
%
%
2x^{10} + x^9 + x^6 + x^4 - x^2 + x + 2 \right)\cr
\tilde{f}_3 =& \textstyle{\frac{1}{3}}(x-1)^2 \left( 
3x^{20}+x^{19}-3x^{18}+x^{17}+6x^{16}-6x^{14}+x^{13}+6x^{12}-x^{11}+\right. \cr
 & \qquad \qquad \qquad \left.- 7x^{10} - x^9 + 6x^8 + x^7 - 6x^6 + 6x^4 + x^3 - 3x^2 + x + 3 \right) \cr
}
$$
}

Neither for $p=2$ nor for $p=3$ the roots of the irreducible factor of
$\tilde{f}_p$ of degree $20$ are integral. Therefore these roots
are not roots of unity and we conclude that $\tilde{f}_p$ has only two
roots that are roots of unity, counted with multiplicities. 
By Corollary \ref{rankNSbound} this implies that the geometric
Picard number of $X_p$ is at most $2$.

Note that besides the hyperplane section $H$, the surface $X_2$ also 
contains the conic $C$ given by $w=g_2=z^2+xy=0$. We have $H^2 = \deg X_2
=4$ and $H \cdot C = \deg C = 2$. As the genus $g(C)$ of $C$ equals
$0$ and the canonical divisor $K$ on
$X_2$ is trivial, the adjunction formula $2g(C)-2 = C\cdot(C+K)$
yields $C^2 = -2$. Thus $H$ and $C$ generate a sublattice
of $\NS(\oX_2)$ with Gram matrix 
$$
\left(
\begin{array}{cc}
4 & 2 \cr
2 & -2 \cr
\end{array}
\right).
$$
\noindent
We conclude that the inner product space $\NS(\oX_2)_{\Q}$ has rank 
$2$ and discriminant $-12 \in \Q^*/{\Q^*}^2$. Similarly, $X_3$
contains the line $L$ given by $w=z=0$, also with genus $0$ and thus
$L^2=-2$. The hyperplane section on $X_3$ and $L$ generate a sublattice
of $\NS(\oX_3)$ of rank $2$ with Gram matrix 
$$
\left(
\begin{array}{cc}
4 & 1 \cr
1 & -2 \cr
\end{array}
\right).
$$
\noindent
We conclude that the inner product space $\NS(\oX_3)_{\Q}$ also has rank 
$2$, and discriminant $-9 \in \Q^*/{\Q^*}^2$.

Let $\rho$ denote the geometric Picard number $\rho = \rk \NS(\oX_h)$. 
It follows from Proposition \ref{nerred} that there is an 
injection $\NS(\oX_h)_{\Q} \hookrightarrow \NS(\oX_p)_{\Q}$ of inner
product spaces for $p=2,3$. Hence we get $\rho \leq 2$. If equality
held, then both these injections would be isomorphisms and
$\NS(\oX_2)_{\Q}$ and $\NS(\oX_3)_{\Q}$ would be isomorphic as 
inner product spaces. 
This is not the case because they have different discriminants. We
conclude $\rho \leq 1$. As a hyperplane section $H$ on $X_h$ has 
self intersection $H^2=4\neq 0$, we find $\rho=1$. 
Since $\NS(\overline{X}_h)$ is a $1$-dimensional even lattice (see
Lemma \ref{NSeven}), the discriminant of $\NS(\overline{X}_h)$ is
even. The sublattice of finite index in $\NS(\oX_h)$ generated by $H$ gives
$$
4= \disc \langle H \rangle = [\NS(\overline{X}_h):\langle H \rangle]^2\cdot
  \disc \NS(\overline{X}_h).
$$
Together with $\disc \NS(\overline{X}_h)$ being even this implies 
$[\NS(\overline{X}_h):\langle H \rangle]=1$, so
$H$ generates $\NS(\overline{X}_h)$, which is isomorphic to $\Pic
\overline{X}_h$ by Lemma \ref{NSeven}.
\end{proof}

\begin{remark}
Corollary \ref{rankNSbound} was pointed out to the author by Jasper Scholten
and people have used it before to bound the
geometric Picard number of a surface. However, since all nonreal roots 
of the characteristic polynomial of Frobenius come in conjugate pairs,
the upper bound has the same parity as the second Betti number of the
surface. For K3 surfaces this means that the upper bound is even (and
therefore at least $2$). The strategy of the proof of Theorem
\ref{picexampletheorem} allows us to sharpen such an upper bound. 
If the reductions modulo two different primes give the same upper
bound $r$, but the corresponding N\'eron-Severi groups have
discriminants that do not differ by a square factor, then in fact
$r-1$ is an upper bound. 

Kloosterman has used our method to construct an elliptic K3
surface with Mordell-Weil rank $15$ over $\overline{\Q}$, see \myref{\kloos}.
In the proof of Theorem \ref{picexampletheorem} we
were able to compute the discriminant up to squares of the 
N\'eron-Severi lattice of $\oX_p$ because we knew a priori a sublattice of 
finite index. Kloosterman realized that it is not always
necessary to know such a sublattice. For an elliptic surface $Y$ over
$\overline{\F_p}$, the image in $\Q^*/{\Q^*}^2$ of the
discriminant of the N\'eron-Severi lattice can also
be deduced from the Artin-Tate conjecture, which has been proved for
ordinary K3 surfaces in characteristic $p\geq 5$, see \myref{\nyog},
Thm.\ 0.2, and \myref{\milneone}, Thm. 6.1.  
It allows one to compute the ratio
$\disc \NS(Y) \cdot \# \Br(Y)/(\NS(Y)_{\tors}^2)$ from 
the characteristic polynomial of Frobenius acting on
$H^2_{\et}(Y,\Q_l)$. For an elliptic surface the Brauer group has 
square order, so this ratio determines the same element in 
$\Q^*/{\Q^*}^2$ as $\disc \NS(Y)$.
\end{remark}

\begin{remark}\label{computerone}
In the proof we counted points over $\F_{p^n}$ for $p=2,3$ and
$n=1,\ldots, 11$ in order to find the traces of powers of Frobenius up
to the $11$-th power. We could have got away with less counting. 
In both cases $p=2$ and $p=3$ we already know a $2$-dimensional
subspace $W$ of $\NS(\oX_p)_{\Q_l}\subset H^2_{\et}(\oX_p,\Q_l)(1)$, generated
by the hyperplane section $H$ and another divisor class.
Therefore it suffices to find out the characteristic
polynomial of Frobenius acting on the quotient $V=H^2_{\et}(\oX_p,\Q_l)(1)/W$.
This implies it suffices to know the traces of powers of Frobenius
acting on $V$ up to the $10$-th power. 

An extra trick was used for $p=3$. The family of planes through the
line $L$ given by $w=z=0$ cuts out a fibration of curves of genus $1$.
We can give all nonsingular fibers
the structure of an elliptic curve by quickly looking for a
point on it. There are efficient algorithms available in for instance
{\sc Magma} to count the number of points on these elliptic curves. 

Using these few speed-ups we let a computer run to
compute the characteristic polynomial of several random surfaces 
given by an equation of the form $wf_1=zf_2$ over $\F_3$ or $wf_1 =
g_1g_2$ over $\F_2$, as in (\ref{pictheeq}). If the middle
coefficient of the characteristic polynomial
was zero, no more effort was spent on trying to find 
the sign of the functional equation (see proof of Theorem
\ref{picexampletheorem}) and the surface was discarded. After one night
two examples over $\F_3$ were found with geometric Picard number $2$
and one example over $\F_2$. With the Chinese
Remainder Theorem this allows us to construct two families
of surfaces with geometric Picard number $1$. One of these families
consists of the surfaces $X_h$. A program written in {\sc Magma}
that checks the characteristic polynomial of Frobenius on
$X_2$ and $X_3$ is electronically available from the author upon request.
\end{remark}

\begin{remark}
For $p=2,3$, let $A_p\subset \NS(\oX_p)$ denote the lattice 
as described in the proof of Theorem \ref{picexampletheorem}, i.e.,
$A_2$ is generated by a hyperplane section and a conic, and $A_3$ is
generated by a hyperplane section and a line. Then in fact $A_p$
equals $\NS(\oX_p)$ for $p=2,3$. Indeed, we have $\disc A_p =
[\NS(\oX_p) : A_p]^2 \cdot \disc \NS(\oX_p)$. For $p=2$ this implies 
$\disc \NS(\oX_2) = -12$ or $\disc \NS(\oX_2) = -3$. The latter is
impossible because modulo $4$ the discriminant of an even lattice of
rank $2$ is congruent to $0$ or $-1$. We conclude $\disc \NS(\oX_2) =
-12$, and therefore $[\NS(\oX_2) : A_2]=1$, so $A_2 = \NS(\oX_2)$.

For $p=3$ we find $\disc \NS(\oX_3) = -9$ or $\disc \NS(\oX_3)=-1$.
Suppose the latter equation held. By the classification of even
unimodular lattices we find that $\NS(\oX_3)$ is isomorphic to
the lattice with Gram matrix 
$$
\left(
\begin{array}{cc}
0 & 1 \cr
1 & 0 \cr
\end{array}
\right).
$$
By a theorem of Van Geemen this is impossible, see \myref{\geem}, 5.4.
From this contradiction we conclude $\disc
\NS(\oX_3) = -9$ and thus $[\NS(\oX_3) : A_3]=1$, so $A_3 =
\NS(\oX_3)$. 
\end{remark}

Since there are 
$\genfrac(){0cm}{1}{4+3}{3}=35$
monomials of degree $4$ in 
$\Q[x,y,z,w]$, the quartic surfaces in $\P^3_\Q$ are parametrized by 
the space $\P^{34}_\Q$, which we will denote by $M$. Let $M'\isom \P^{27}
\subset M$ denote the subvariety of those surfaces $X$ for which the
coefficients of the monomials $x^4$, $x^3y$, $x^3z$, $y^4$, 
$y^3x$, $y^3z$, and $x^2z^2$ in the defining polynomial of $X$ are all
zero. Note that the vanishing of the coefficients of the first six of these
monomials is equivalent to the tangency of the plane $H_w$ given by $w=0$ to
the surface $X$ at the points $P=[1:0:0:0]$ and $Q=[0:1:0:0]$. Thus,
the vanishing of these coefficients
yields a singularity at $P$ and $Q$ in the plane curve $C_X = H_w \cap
X$. If the singularity at $P$ in $C_X$ is not worse than a double point, 
then the vanishing of the coefficient of $x^2z^2$ is equivalent to
the fact that the line given by $y=w=0$ is one of the limit-tangent lines to
$C_X$ at $P$. 

\begin{proposition}\label{picopenU}
There is a nonempty Zariski open subset $U \subset M'$ such that every
surface $X \in U$ defined over $\Q$ is smooth and 
has infinitely many rational points. 
\end{proposition}
\begin{proof}
The singular surfaces in $M'$ form a closed subset of $M'$. So
do the surfaces $X$ for which the intersection $H_w \cap X$ has worse
singularities than just two double points at $P$ and $Q$. Leaving out
these closed subsets we obtain an open subset $V$ of $M'$. Let $X \in
V$ be given. The plane quartic curve $C_X = X \cap H_w$ has two double
points, so the geometric genus $g$ of the normalization $\tilde{C}_X$
of $C_X$ equals $p_a-2$, where $p_a$ is the arithmetic genus of $C_X$,
see \myref{\hag}, exc. IV.1.8. As we
have $p_a = \frac{1}{2}(4-1)(4-2)=3$, we get $g=1$. 
Now assume $X$ is defined over $\Q$. 
One of the limit-tangents to $C_X$ at $P$ is given by $w=y=0$. Its
slope, being rational, corresponds to a rational point $P'$ on
$\tilde{C}_X$ above $P$. Fixing  
this point as the unit element $\O = P'$, the curve $\tilde{C}_X$
obtains the structure of an elliptic 
curve. Let $D \in \Pic^0(\tilde{C}_X)$ be the pull back under
normalization of the divisor $P-Q \in \Pic^0(C_X)$. By the theory of
elliptic curves there is a unique point $T$ on $\tilde{C}_X$ such that 
$D$ is linearly equivalent to $T-\O$, see \myref{\silv}, Prop.\ III.3.4.
As $D$ is defined over $\Q$, so is $T$. By Mazur's
theorem (see \myref{\silv}, Thm. III.7.5 for statement, 
\myref{\maz}, Thm. 8 for a proof),
the point $T$ has finite order if and only if
$mT=\O$ for some $m \in \{1,2,\ldots, 10, 12\}$. Note that we have 
$\lcm(1,2,\ldots,10,12)=2520$. Take for $U$ the
complement in $V$ of the closed subset of those $X$ for which we have
$2520T = \O$ for the corresponding point $T$ on $\tilde{C}_X$. Then
each $X \in U$ contains an elliptic curve with infinitely many
rational points.
By choosing a Weierstrass equation, one verifies easily that if we
take $X_0 = X_h$ with $h=0$, then the corresponding point $T$ on
$\tilde{C}_{X_0}$ satisfies $mT \neq \O$ for $m \in \{1,2,\ldots, 10, 12\}$.
Therefore, we find $X_0 \in U$, so $U$ is nonempty.
\end{proof}

\begin{remark}
If $\tilde{C}_X$ is the normalization of $C_X$ as in the proof
of Proposition \ref{picopenU}, then generically there is another rational
point $P''$ on $\tilde{C}_X$ above $P$, besides $P'$. Generically this
point also has infinite order and the Mordell-Weil rank of
$\tilde{C}_X$ is at least $2$ with independent points $P''$ and $T$ as
in the proof of Proposition \ref{picopenU}. 
For $X=X_h$ with $h=0$ however, the curve $\tilde{C}_X$ is given by 
$$
3x^2y^2 + xy^2z + 4xyz^2 + 2xz^3 + 5yz^3 + z^4 = 0.
$$
As the point $P=[1:0:0]$ is a cusp, there is only one point above $P$
on $\tilde{C}_X$ here. The conductor of this elliptic curve
equals $686004$. Both points on $\tilde{C}_X$ above
$Q=[0:1:0]$ are rational and we have an extra rational point $[1:1:-1]$.
These generate the full Mordell-Weil group of rank $3$.
\end{remark}

\begin{remark}
By requiring other coefficients to vanish than is required for $M'$, we
can find quartic surfaces $Y$ for which the plane $H_w$ given by
$w=0$ is tangent at $[1:0:0:0]$, $[0:1:0:0]$, and $[0:0:1:0]$. Then 
the intersection $H_w \cap Y$ has geometric genus $0$ and if its
normalization has a point defined over $\Q$, then 
this intersection is birational to $\P^1$. The quartic surface $Z$ given by
\begin{equation}\label{picexampleP1}
w(x^3+y^3+z^3+x^2z+xw^2)=3x^2y^2-4x^2yz+x^2z^2+xy^2z+xyz^2-y^2z^2
\end{equation}
is an example of such a surface.
As in the proof of Theorem \ref{picexampletheorem}, modulo $3$ the 
surface $Z$ contains
the line $z=w=0$. Also, the reduction of $Z$ at $p=2$ contains a conic
again, as the right-hand side of 
(\ref{picexampleP1}) factors over $\F_4$ as $(xy+xz+\zeta
yz)(xy+xz+\zeta^2 yz)$, with $\zeta^2+\zeta+1=0$. An argument very
similar to the one in the proof of Theorem \ref{picexampletheorem}
then shows
that $Z$ also has geometric Picard number $1$ with the Picard group
generated by a hyperplane section. The only difference is that 
Frobenius does not act trivially on the conic $w=xy+xz+\zeta yz =0$. 
The hyperplane section $H_w \cap Z$ is a curve of geometric genus $0$,
parametrized by 
$$
[x:y:z:w] = [-(t^2+t-1)(t^2-t-3): 2(t+2)(t^2+t-1): 2(t+2)(t^2-t-3):0].
$$
The Cremona transformation $[x:y:z:w] \mapsto [yz:xz:xy]$ gives a
birational map from this curve to a nonsingular plane curve of degree $2$. 
It turns out that the curve on $Z$ given by $x=0$ has a triple point
at $[0:0:0:1]$, so it is birational to $\P^1$ as well. It can be parametrized
by 
$$
[x:y:z:w] = [0:1+t^3:t(1+t^3):-t^2].
$$
\end{remark}

From the local and global Torelli theorem for K3 surfaces, see
\myref{\pjatshaf}, one can find a very precise description of the
moduli space of polarized K3 surfaces in general, see
\myref{\beau}. \check{find better references}
A polarization of a K3 surface $Z$ by a very ample divisor of degree $4$
gives an embedding of $Z$ as a smooth quartic surface in $\P^3$ with
the very ample divisor corresponding to a hyperplane section. 
An isomorphism between two smooth quartic surfaces in $\P^3$ that
sends one hyperplane section to an other hyperplane section comes from an
automorphism of $\P^3$. As any two hyperplane sections are linearly
equivalent, we conclude that
the moduli space of K3 surfaces polarized by a very ample divisor of
degree $4$ is isomorphic to the open subset in $M=\P^{34}$ of smooth 
quartic surfaces modulo the
action of $\PGL(4)$ by linear transformations of $\P^3$.
We are now ready to prove the main theorem of this article.
\smallskip

\begin{proof}[Proof Theorem \ref{picmain}]
%
By the description of the moduli space of K3 surfaces polarized by a
very ample divisor of degree $4$ given above,
%
%
it suffices to prove that the set $S\subset M(\Q)$ of smooth surfaces 
with geometric Picard number $1$ and infinitely many rational points
is Zariski dense in $M$. 
%
%
%
%
We will first show that $S \cap M'$ is dense in $M'$. Note that the
coefficients of the monomials $x^4$, $x^3y$, $x^3z$, $y^4$, $y^3x$,
$y^3z$, and $x^2z^2$ in $wf_1+2zf_2-3g_1g_2$ in (\ref{pictheeq}) are
zero, so if the coefficients of these monomials in a homogeneous
polynomial $h \in R$ of degree $4$ are all zero, then $X_h$ is
contained in $M'$. It follows that the set 
$$
T = M' \cap \{X_h : h \in R, h \mbox{ homogeneous of degree } 4\}
$$
is dense in $M'$. Let $U$ be as in Proposition \ref{picopenU}. Then
$U$ is a dense open subset of $M'$, so $T \cap U$ is also dense in
$M'$. By Theorem \ref{picexampletheorem} and  Proposition
\ref{picopenU} every surface in $T \cap U$ has geometric Picard number
$1$ and infinitely many rational points. Thus we have an inclusion $T
\cap U \subset S \cap M'$, so $S \cap M'$ is dense in $M'$ as well.

Let $W$ denote the vector space of $4\times 4$--matrices over $\Q$ and let 
$T$ denote the dense open subset of $\P(W)$ corresponding to elements
of $\PGL(4)$. Let $\varphi \colon T \times M' \rightarrow M$ be
given by sending $(A,X)$ to $A(X)$. Note that $T(\Q) \times (S \cap M')$ is 
dense in $T \times M'$ and $\varphi$ sends $T(\Q) \times (S \cap M')$ to
$S$. Hence, in order to prove that $S$ is dense in $M$, it suffices to
show that  
$\varphi$ is dominant, which can be checked after extending to the
algebraic closure. A general quartic surface in $\P^3$ has a
one-dimensional family of bitangent planes, i.e., planes that are
tangent at two different points. This is closely related to the theorem
of Bogomolov and Mumford, see the appendix to \myref{\morimukai}.
\check{reference}
In fact, for a general quartic surface $Y \subset \P^3$, there is such
a bitangent plane $H$, for which the two tangent points are ordinary
double points in the intersection $H \cap Y$. Let $Y$ be such a
quartic surface and $H$ such a plane, say tangent at $P$ and
$Q$. Then there is a linear transformation that sends $H$, $P$, and
$Q$ to the plane given by $w=0$, and the points $[1:0:0:0]$ and
$[0:1:0:0]$. Also, one of the limit-tangent lines to the curve $Y \cap H$ at
the singular point $P$ can be sent to the line given by $y=w=0$. This
means that there is a linear transformation $B$ that sends $Y$ to an
element $X$ in $M'$. Then $\varphi(B^{-1},X)=Y$, so $\varphi$ is
indeed dominant.
\end{proof}

\begin{remark}
The explicit polynomials $f_1,f_2,g_1,g_2$ for $X_h$ in (\ref{pictheeq}) were
found by letting a computer pick random polynomials modulo $p=2$ and $p=3$
such that the surface $X_h$ with $h=0$ is contained in $M'$ as in Proposition 
\ref{picopenU}. The computer then computed the characteristic polynomial
of Frobenius and tested if there were only $2$ eigenvalues that were
roots of unity, see Remark \ref{computerone}. 
\end{remark}

\begin{remark}
In finding the explicit surfaces $X_h$ not much computing power was
needed, as we constructed the surface to have good reduction at small 
primes $p$ so that counting points over $\F_{p^n}$ was relatively easy. 
Based on ideas of for instance Alan Lauder, Daqing Wan, Kiran
Kedlaya, and Bas Edixhoven, it should be possible to develop more
efficient algorithms for finding characteristic polynomials of (K3)
surfaces.  \check{first names? reference?} Together with
these algorithms, the method used in the proof of Theorem
\ref{picexampletheorem} becomes a strong tool in finding Picard numbers
of K3 surfaces over number fields. 
\end{remark}

\section{Open problems}

We end with the remark that still very little is known about the
arithmetic of K3 surfaces, especially those with
geometric Picard number $1$. We reiterate three questions that remain
unsolved. 

\begin{question}
Does there exist a K3 surface over a number field such that the set of
rational points is neither empty nor dense?
\end{question}

\begin{question}
Does there exist a K3 surface over a number field with geometric
Picard number $1$, such that the set of rational points is potentially dense?
\end{question}

\begin{question}
Does there exist a K3 surface over a number field with geometric
Picard number $1$, such that the set of rational points is not
potentially dense? 
\end{question}

\begin{acknowledgements}
The author thanks the American Institute of Mathematics (Palo
Alto) and the Institut Henri Poincar\'e (Paris) for inspiring working
conditions. The author also thanks Bjorn Poonen, Arthur Ogus, Jasper
Scholten, Bert van Geemen, and Hendrik Lenstra
for very useful discussions, and Brendan Hassett for pointing out a
mistake in the first version of this article.
\end{acknowledgements}

\end{document}